# The Calculus Concept Readiness (CCR) Instrument:
## Assessing Student Readiness for Calculus

Marilyn Carlson, Arizona State University

Bernard Madison, University of Arkansas

Richard West, Francis Marion University

**Abstract**

The Calculus Concept Readiness (CCR) instrument is based on the broad body of mathematics education research that has revealed major understandings, representational abilities, and reasoning abilities students need to construct in precalculus level courses to be successful in calculus. The CCR is a 25-item multiple-choice instrument, and the CCR taxonomy articulates what the CCR assesses. The methodology used to develop and validate the CCR is described and illustrated. Results from administering the CCR as a readiness examination in calculus are provided along with data to guide others in using the CCR as a readiness examination for beginning calculus.

**Introduction**

A commonly cited reason for high attrition from precalculus to calculus is that standards for what should be taught in courses like precalculus mathematics are unclear and inconsistent (Arismendi-Pardi, 1996; Schmidt, Houang, & Cogan, 2002). There is now substantial research on what is involved in learning key ideas of algebra through beginning calculus (e.g., Dubinsky & Harel, 1992; Thompson, 1994a, 1994b, 1994c; White & Mitchelmore, 1996; Carlson & Rasmussen, 2008; Thompson, 2008). However, a cursory examination of the commonly used curricula suggests that this research has had little influence on what is being taught, resulting in many students enrolling in beginning calculus without the prerequisite knowledge. This manuscript provides a response by articulating key reasoning abilities, understandings, and representational abilities that are critical for learning calculus. It also presents the Calculus Concept Readiness (CCR) instrument that can be used to assess both student readiness for calculus and the effectiveness of specific precalculus curriculum or courses in preparing students to learn key calculus concepts.

*Foundational Understandings for Learning Calculus*

National organizations and policy makers have made repeated calls for precalculus level curricula to place greater emphasis on functions (College Entrance Examination Board, 1959; Hamley, 1934; Hedrick, 1922;  Klein, 1883; National Council of Teachers of Mathematics, 1934, 1989, 2000). Over the past 25 years many mathematics education researchers have found that student difficulty in understanding key ideas of calculus are rooted in their weak understanding of the function concept (Tall & Vinner, 1981; Williams, 1991; Dubinsky et al., 1992; Kaput, 1992; Monk, 1992; Tall, 1992, 1996; Breidenbach, Dubinsky, Hawks, Nichols, 1991; Thompson, 1994a; Cottrill et al., 1996; Carlson, 1998; Zandieh, 2000; Oehrtman, 2004, 2008a, 2008b; Engelke, 2007; Smith, 2008) and their inability to use functions to reason about and represent the relationships between quantities and how they change together (Monk & Nemirovksy, 1994; Thompson, 1994b; Carlson et al., 2002). Other studies have revealed that calculus students' narrow conception of variable and quantity (Jacobs, 2002, Thompson, 1994; Ursini & Trigueros, 2001), and their inability to reason about and represent how two quantities change together cause problems for students in learning ideas of calculus such as limit, derivative, accumulation and the Fundamental Theorem of Calculus (Monk, 1987; Kaput, 1992; Thompson, 1994c; Cottrill, Dubinsky, Nichols, Schwingendorf, Thomas, Vidakovic, 1996; Carlson et al., 2002; Carlson, Smith, Persson, 2003; Engelke, 2007; Oehrtman, Carlson & Thompson, 2008; Smith, 2008) and ideas in differential equations (Rasmussen, 2000).

Early studies of student understandings of the function concept revealed common misconceptions that students encounter when learning ideas of precalculus and beginning calculus (e.g., Monk, 1992; Sierpinska, 1992; Vinner & Dreyfus, 1989), including students' narrow conception of variable as something to solve for and their strong tendency to view a graph as a picture of an event, rather than a representation of how two quantities change together. Investigations of precalculus students during and after instruction have provided insights into the reasoning abilities and understandings that support student development of an understanding of ideas of function (Dubinsky et al., 1992; Carlson, 1998), function composition (Carlson, 1998; Engelke, Oehrtman & Carlson, 2005; Engelke, 2007), function inverse (Dubinsky et al., 1992; Carlson, 1998), quantity (Moore, Carlson & Oehrtman, submitted), exponential growth (Strom, 2008; Castillo-Garsow, 2010), and central ideas of trigonometry (Moore, 2010). These studies have consistently revealed that student success in learning calculus requires that students view a function, in each of its representational forms (words, formula, graph, table), as a means of associating input values from a function's domain to output values in the function's range for the purpose of representing how values of one quantity change with the values of another quantity (e.g., Tall, 1992; Thompson, 1994a; Cottrill et al., 1996; Carlson et al., 2002; Engelke, 2007; Carlson & Ramussen, 2008; Oehrtman, 2008a, 2008b; Smith, 2009; Moore, 2010).

This body of literature informed the development of the Calculus Concept Readiness (CCR) taxonomy. The taxonomy then informed the work of a committee of mathematicians and mathematics educators in developing or selecting CCR items[1].

## The Process of Developing the CCR Instrument

The development of the CCR was inspired and informed by mathematics education research and instrument validation methods of the Precalculus Concept Assessment (PCA) (Carlson, Oehrtman & Engelke, 2010) and by validation methods used in developing the Mechanics Baseline Test (Hestenes & Wells, 1992) and the Force Concept Inventory (Hestenes, Wells, & Swackhammer, 1992). The development of these instruments began with the creation of initial taxonomies to characterize broad categories of student thinking related to concepts that had been identified by research, and mathematicians or scientists to be important. These constructs were gleaned from mathematics education research literature and were described and articulated in the instrument taxonomies[2]. Open-ended questions were designed to assess these reasoning abilities and understandings. Question wording and item distractors for the multiple-choice versions were developed based on student interview data that illustrated common student thinking when responding to each open-ended item. Further interviews were conducted with students after the initial multiple-choice items were developed to identify or refine item distractors so that they were representative of the five most common student responses, and to verify that questions and answer choices were interpreted correctly. The taxonomy and CCR items went through multiple cycles of refinement. Clinical interviews with students were conducted repeatedly until each CCR item had been validated to: i) be consistently interpreted, ii) assess the knowledge intended by the item designer, and iii) have optimal distractors that were representative of the student thinking revealed during the interviews.

---

[1] The CCR instrument is part of the Placement Testing Suite of the Mathematical Association of America that is delivered by Maplesoft. Currently, parallel forms of the initial CCR are being developed, as are samples of items for national dissemination to inform precalculus instruction.

*What Does CCR Assess?*

The CCR taxonomy (Table 1) articulates key reasoning abilities and understandings assessed by the 25-item multiple-choice Calculus Concept Readiness (CCR) instrument. Eighteen of the twenty-five CCR items assess or rely on student understanding of the concept of function. Five items assess student understanding or use of trigonometric functions, and four items assess student understanding or ability to use exponential functions. While ten items are situated in an applied (or word problem) context and require students to reason about quantities and use ideas of function, function composition, or function inverse to represent how the quantities change together.  The CCR taxonomy is organized into categories of reasoning abilities, conceptual understandings, representing and interpreting growth patterns in two quantities, ideas central to trigonometry, and other general abilities that have been found to be foundational for learning calculus.

TABLE 1
**The CCR Taxonomy**

Reasoning abilities
   R1    Use proportional reasoning
   R2    View a function as a process that defines how input values in a function's domain
            are mapped to or associated with output values in a function's range (process
            view of function)

---

[2] According to Lissitz and Samuelsen (2007), the development of a valid examination should always begin by identifying the constructs worthy of assessment.

R3    Reason about and represent how two quantities change together (engage in covariational reasoning)

Understand and use the following concepts or ideas

U1    Quantity
U2    Variable
U3    Slope/Constant rate of change
U4    Average rate of change
U5    Function composition
U6    Function inverse
U7    Function translations (horizontal and vertical shifts)

Understand, represent and interpret function growth patterns

G1    Linear
G2    Exponential
G3    Non-linear polynomial
G4    Rational
G5    Periodic

Understand central ideas of trigonometry

T1    Angle measure
T2    Radian as a unit of measure
T3    Sine and cosine functions as the covariation of an arc's length (measured in units of the circle's radius) and the horizontal or vertical coordinate of the arc's terminus (measured in units of the arc's radius)[3]
T4    Sine and cosine functions as a representation of the relationship between an angle measure and sides of a right triangle

Other abilities

A1    Solve equations
A2    Represent and interpret inequalities
A3    Represent and interpret absolute value inequalities
A4    Use and solve systems of equations
A5    Understand and use function notation to express one quantity in terms of another

---

[3] These questions exploit the idea that every circle can be considered a unit circle.

**Calculus Concept Readiness Items**

The following section presents three CCR items and discusses the primary understandings and reasoning abilities needed to provide a correct response. Select data from validating the answer choices for these items are also provided along with quantitative data from administering the items to 631 students.

*A Function Composition Item*

  A function composition word problem (Table 2) prompted students to define the area of a circle in terms of its circumference. Interviews conducted with students who answered the question in an open-ended format revealed that the item assesses students' ability to identify the quantities to be related (U1) in a word problem. The item also assesses student ability to interpret the phrase, 'express area $A$ as a function of circumference $C$' as a prompt to construct a formula of the form, $A = $ <some expression that contains C> (A5). Students must recall the formulas for the area and circumference of a circle, and view these formulas/functions as a means of determining values for one quantity when the values of another are known (R2). To obtain the function that relates the area and circumference of a circle, students must recognize that these formulas can be combined by inverting one formula (U5), then composing (U4) the two function formulas for the purpose of relating the quantities area and circumference.

TABLE 2

Area-Circumference Item

Which of the following formulas defines the area, $A$, of a circle as a function of its circumference, $C$?

a.   $A = \dfrac{C^2}{4\pi}$

b.   $A = \dfrac{C^2}{2}$

c.   $A = (2\pi r)^2$

d.   $A = \pi r^2$

e.   $A = \pi(\dfrac{1}{4}C^2)$

The answer choices for the Area-Circumference item were initially identified by administering the item in an open ended format on a written examination. Follow-up interviews were conducted with 19 students, 7 students who provided the correct answer and 3 students who constructed each of the four most common incorrect responses. The students were prompted to explain their approach and the thinking they used to construct their solutions.

Analysis of interview data with students who provided a correct answer revealed a common approach of verbalizing the problem goal (to express area in terms of circumference). These students typically wrote the formulas $A = \pi r^2$ and $C = 2\pi r$, and eventually recognized the need to re-express the circumference formula by solving for $r$ (or determining the inverse of the circumference formula) to obtain $r = \dfrac{C}{2\pi}$. When students were prompted to explain how they knew to solve the circumference formula for $r$, a common response as articulated by one student was, "Since $A = \pi r^2$ and I need a formula for area in terms of circumference, I can solve $C = 2\pi r$ for $r$. So now I can put $r = \dfrac{C}{2\pi}$ in for $r$, since this is the same as $r$. This will give me a function that takes values of the circumference $C$ and computes the area $A$." Conceptualizing an approach for linking two function processes for the purpose of expressing one quantity in terms of another is a key construction that influences whether students are able to represent, using a formula or graph, the relationship between quantities in related rate problems in beginning calculus (Carlson, 1998; Engelke, 2007). In Engelke's study, students sometime referred to function composition as determining the middle variable that links two other variables. The students expressed that their goal was to re-express the variable that linked the two variables to be related in such a way that values of the input variable (the first quantity to be related) passed through the middle variable to produce a value for the output variable (the second quantity to be related).

Students who constructed the answers stated in response choices (b), (c), and (d) initially appeared to be making careless errors. However, their explanations during the interviews revealed that their careless errors were typically rooted in a procedural approach to responding to the problem. They did not view the substitution of $\dfrac{C}{2\pi}$ for $r$ as the linking of two function processes, nor did they view $\dfrac{C}{2\pi}$ as a quantity that represented the value of $r$ (given values for $C$) that needed to be squared. Students who chose answer (d) merely wrote down the formula for the area of a circle in terms of the radius, and did not recognize that a request to construct a formula to represent area $A$ of a circle in terms of its circumference $C$ implied that they needed to construct a formula with $A$ equal to some expression written in terms of $C$. These students were clearly not thinking about function composition as a means of stringing together two function processes.

When administering the validated version of this multiple choice item to 631 beginning calculus students at two large public universities and one large private university during the first week of classes, only 28% of these students selected the correct response. Since at least 13% of the students (out of 631) selected each answer choice, this quantitative data (Table 3) further supports the findings from analyzing the qualitative data that the five answer choices are representative of common student answers.

TABLE 3

Area-Circumference Answer Percentage

| a. 28.30% | b.13.67% | c. 20.19% | d. 21.62% | e. 15.58% |
|-----------|----------|-----------|-----------|-----------|

*A Quantitative Reasoning and Slope Item*

The following item requires that students conceptualize two quantities (the distance of the top of the ladder from the floor and the distance of the base of the ladder from the wall) (U1). They then need to imagine how the distance of the base of the ladder from the wall changes as the top of the ladder increases to twice its original distance from the floor (R3). As they imagine how these measurements change together (engage in covariational reasoning), they also need to think about how the ratio of the changes in these two quantities (slope) (U3) changes as the distance of the ladder from the wall decreases (R3). Another correct justification involved imagining how the ratio that represents *the slope of a line* changes if the algebraic form of the slope is changed by doubling the value of the numerator and decreasing the value of the denominator. Students who used this approach verbalized that doubling the numerator would produce a slope twice what it was before, but doubling the numerator and decreasing the value of the denominator would produce a slope that is more than twice as large as what it was before. These students were able to imagine how concurrent changes in two lengths on the illustration impacted the values of the numerator and denominator and how these values impacted the value of the ratio.

TABLE 4

Slope of Ladder Item

A ladder that is leaning against a wall is adjusted so that the distance of the top of the ladder from the floor is twice as high as it was before it was adjusted.

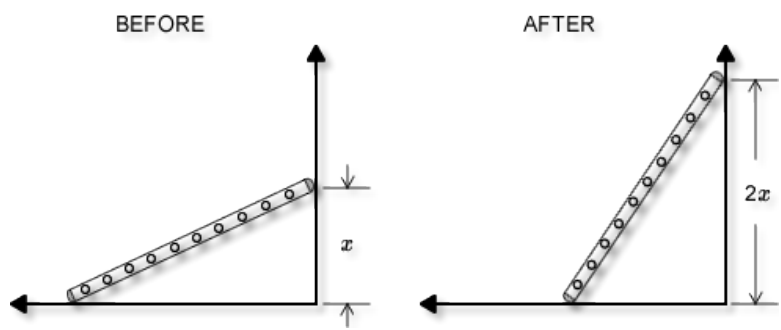

The slope of the adjusted ladder is:

a.   Less than twice what it was
b.   Exactly twice what it was
c.   More than twice what it was
d.   The same as what it was before
e.   There is not enough information to determine if any of a through d is correct.

TABLE 5

Slope of Ladder Answer Percentages

a) 5.89%          b) 48.01%          c) 27.34%          d) 3.50%          e) 14.63%

Only 27% of 631 beginning calculus students provided a correct response to the slope item.  The most common incorrect response was choice (b), exactly twice what it was. Interviews with 7 students who selected this choice revealed that these students were not conceptualizing the slope of the straight line (ladder) as a ratio of two quantities. The interviews revealed that the students did not consider how the ratio of *the increase in the top of the ladder by a factor of 2* and the decrease in the *distance of the base of the ladder from the floor* affected the ladder's slope. Rather, they replied with responses that indicated they were only focusing on the amount of increase of the top of the ladder.  When prompted to elaborate on their rationale for choice (b), one student explained, "since the top of the ladder is two times higher when the ladder is repositioned, the slope will be twice as large."  The student failed to consider the effect of the shortened distance of the base of the ladder from the wall and never considered the slope as representing a ratio of the two distances. Students who selected answer (e), not enough information to determine, thought they needed exact numbers for the distances to be able to compute and compare the slopes. Students who thought the slope was less than twice what it was (answer (a)) thought that the smaller denominator made the slope less than twice what it was. Students who indicated that the ladder's slope had not changed (answer (d)) indicated that the ladder itself did not change because the positioning of the ladder on the wall had no effect on the shape of the ladder. These students were not thinking about the slope as representing a ratio of the change in the output variable and the change in the independent variable.

*A Trigonometry Item: Assessing Ideas of Angle Measure and Sine Function*

The following trigonometry item (Table 6) relies on students having developed robust conceptions of function (R2), angle measure (T1), radian (T2) and the sine function (T3), in addition to being able to imagine and fluently reason about how two quantities change in tandem (an angle measure and distance) (R3). Responding to this item requires that students first have a conception of the sine function as representing the covariation of an angle measure and the distance of a point on the unit circle (positioned at the end of the terminal side of the angle) from a horizontal line through the circle's center (T3, R3). Determining the function that expresses $d$ in terms of $k$ will also require that students conceptualize an angle measure as an arc length that is cut off by the rays of an angle positioned at the center of a circle and measured in lengths of the circle's radius (T1 and T3) ($k$ needs to be divided by 47). As students imagine $k$ varying, they also imagine how $d$ varies (R3) and recognize that, since the output of the sine function is measured in lengths of the radius (radian), this value needs to be multiplied by 47 to express the value of $d$ in feet.

TABLE 6

Periodic Motion Item

Starting at P and ending at Q, an object travels counterclockwise $k$ feet along a circle with radius 47 feet. If $d$ represents the directed distance (in feet) from the horizontal diameter to Q, which of the following could express $d$ as a function of $k$?

a.   $d = f(k) = 47 \sin (k)$

b.   $d = f(k) = 47 \sin (47k)$

c.   $d = f(k) = 47 \sin (\dfrac{k}{94})$

d.   $d = f(k) = 47 \sin (\dfrac{47}{k})$

e.   $d = f(k) = 47 \sin (\dfrac{k}{47})$

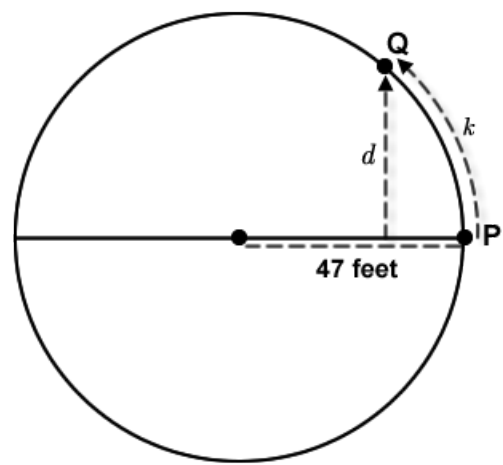

Only 21% of 631 beginning calculus students selected the correct answer, (e) (TABLE 7). As noted above, responding to this question requires that students understand ideas of angle measure, radian and sine function, and can reason about quantities, their variation and their covariation (i.e., how the values of two quantities, $k$ and $d$, change in tandem).

TABLE 7

Periodic Motion Answer Percentages

| a) 24.96% | b) 8.27% | c) 19.40% | d) 20.19% | e) 21.14% |
|---|---|---|---|---|

Students who selected answer (a) did not understand the idea of angle measure and that an angle measure can be expressed in lengths of the radius. The interview data also revealed that these students did not understand that the input to the sine function must be expressed as an angle measure. Students who selected answer (b) did not understand the idea of angle measure. Students who chose answer (c), (d), or (e) also had a weak understanding of angle measure and radian as a measure of the number of radius lengths subtended by the rays of an angle positioned at a circle's center.

## Item Analysis of Calculus 1 Results

During the fall semester of 2009, the CCR was administered to 215 Calculus 1 students at a large public university during the first week of class.  Most of these students were freshmen and were placed in the course by virtue of their ACT mathematics test scores. The data generated from this pilot (Table 8) was used to determine the correlation of CCR scores with other measures of student preparedness and performance in calculus. The data reveal that a student's CCR score is most highly correlated with her ACT mathematics score. At the same time reasonably strong correlations exist between CCR scores and course grades and test scores.

TABLE 8

Correlation of CCR Score with Other Measures

| Measure | Number of Records | Correlation |
|---|---|---|
| ACT Mathematics Score | 144 | 0.54 |
| First Test Score | 207 | 0.50 |
| Mid-Term Test Score | 195 | 0.41 |
| Final Course Grade | 215 | 0.51 |

*Validity Measures*

As noted above, analysis of CCR scores from assessing 215 Calculus I students at the beginning of the course shows a reasonably strong connection between levels of CCR scores and success in their Calculus I course, as measured by course grades. The mean CCR score of students whose course grades were A, B, or C is 11.83, while the mean CCR score for students whose course grades were D, W, or F is 8.49.  The table below (Table 9) shows that CCR scores of 0-7 correlate with a success rate of about 30%, CCR scores of 8, 9, or 10 correlate with a success rate of about 50%, and CCR scores of 11 or more correlate with a success rate in excess of 66%.  It is also noteworthy that 100% of students who received a score of 17 or greater received a course grade of A, B, or C. This finding suggests that CCR items assess prerequisite knowledge for calculus.

TABLE 9

| CCR Score | N | % A, B, or, C | % D, W, or F |
|-----------|-----|---------------|---------------|
| 0 - 7 | 23 | 30 | 70 |
| 8 | 15 | 44 | 56 |
| 9 | 16 | 50 | 50 |
| 10 | 16 | 47 | 53 |
| 11 - 16 | 125 | 66 | 34 |
| 17 - 25 | 15 | 100 | 0 |

Other pilots of the CCR at public and private universities show similar results. Some of those are summarized below along with cautionary advice on interpreting placement test results. More data will be available in the upcoming year to guide institutions in determining cut-scores for requiring or advising students to complete a prerequisite precalculus course[4].

The pilot testing described above illustrates one approach for validating a placement instrument. There are limitations in attempting to correlate a student's CCR score with her course grade since factors other than knowledge of and facility with precalculus concepts influence grades. As a result, it is unlikely that a test administered at the beginning of a course will be strongly correlated with a student's course grade. Secondly, students are usually placed in Calculus 1 by some criteria. In the case above, ACT mathematics scores were used along with a back-up placement test score when the ACT test score was missing or flawed. These conditions likely reduced the variation in the CCR scores of this population, resulting in the likelihood that the predictive power of CCR scores was reduced.

Another type of validity study involved comparison of the CCR results to the outcome of the prerequisite precalculus course. Since students who earn grades of A, B, or C in a precalculus course are presumably ready to learn the content taught in a Calculus 1 course, testing precalculus students at the end of the course should reveal a correspondence between success in the course and a CCR score. Data from 31 precalculus students at a large public university were tested at the end of their course. The correlation between CCR test scores and course grades in precalculus was 0.58 with the following results that, when combined with the predictive study, yields reasonably consistent results.

---

[4] Previous studies (Carlson, Oehrtman & Engelke, 2010) suggest that a precalculus course that focuses primarily on procedures and skills may not lead to shifts in students' scores on a conceptually oriented readiness examination.

TABLE 10

| CCR Score | Number of Scores | % A, B, or C |
|-----------|------------------|--------------|
| 0 - 7 | 6 | 33 |
| 8 - 10 | 14 | 57 |
| 11 - 25 | 11 | 91 |

At this university, as at most universities, students who are "on track" enroll in Calculus 1 in the fall of their freshman year.  Calculus 1 in the second semester has a student population with a larger variety of backgrounds: some were unsuccessful in Calculus 1 in the fall, some completed university or college precalculus courses, and some are entering freshmen, after a semester or more out of school. Classes of Calculus 1 in the spring are generally weaker, and the CCR results from testing 281 Calculus 1 students in the spring of 2010 revealed lower correlations than for fall semester students. The mean scores were less (8.71 versus 10.38) and grades were lower (1.36 out of 4 versus 1.79 out of 4).  Reliability estimates were lower (0.54 versus 0.66) and correlations with final grades were lower (0.28 versus 0.51), implying that placement tests are not as effective in placing spring semester students. Other factors such as motivation and study habits may be playing a greater role in student success during spring semesters, also contributing to the generally lower success rates. It is also likely that CCR is more effective in predicting success in Calculus 1 courses that include a focus on understanding and using concepts of calculus[5]. The CCR scores and grades (of A, B, or C) for these spring semester students are summarized (Table 11).

Table 11

| CCR Score | Number of Scores | % of A, B, or C |
|-----------|------------------|-----------------|
| 0 – 7 | 98 | 38 |
| 8 – 10 | 81 | 54 |
| 11 – 25 | 71 | 63 |

*CCR Testing at Two Other Universities*

In addition to the pilot testing of CCR at one large public university, CCR was also pilot-tested at another large public university and at a large private university.  The institutional settings of the three universities were different, and the results reflect some of those differences.

_______________________

[5] Recall that CCR items assess student understandings of precalculus ideas that are documented to be foundational for understanding and using key ideas of calculus.

The private university has a more restrictive standard for admitting students to the university. It also places students in Calculus 1 using scores on a departmental placement test with some weight given to the ACT or SAT score. As expected, the CCR scores were higher, with a 1-point–higher mean score as compared to the two public universities. The data also revealed that student grades were generally higher.

TABLE 12

| CCR Score | Number of Scores | % A, B, or C |
|-----------|-----------------|--------------|
| 0 - 7     | 11              | 40           |
| 8 - 10    | 21              | 81           |
| 11 - 25   | 36              | 86           |

At the other large public university in our study, 57 students in three Calculus 1 sections were tested using CCR. Students were placed in the course using a placement test score. The sections in which the participating students were enrolled also included collaborative learning and extra support sessions that may have accommodated for student weaknesses and influenced student success. Of the 18 students with CCR scores of 0 – 9, 66% earned course grades of A, B, or C, and 95% of the 39 students with CCR scores of 10 – 25 earned grades of A, B, or C.

These results revealed that students who receive higher CCR scores generally performed better in Calculus 1. In fact, students who scored 11 or higher on CCR passed calculus at rates of 66%, 86% and 91% at the three universities reported in this study. This finding is consistent with the results reported by Carlson et al. (2010) using the Precalculus Concept Assessment (PCA). The authors reported that 77% of 248 Calculus I students who took PCA at the beginning of a fall semester, and received a score of 13 (out of 25) or higher, were awarded a course grade of A, B, or C. However, the 0.51 correlation coefficient between the initial CCR score and final course grade, was slightly higher than the correlation coefficient of 0.47 between PCA scores, also administered at the beginning of a fall semester, and final course grade. The success rate of 95% or greater for students who scored a 17 or higher on CCR suggests that collectively the items on CCR assess essential understandings that are used in beginning calculus.

## Concluding Remarks

The CCR examination provides a broad snapshot of student understanding of central ideas of precalculus as articulated in the CCR Taxonomy. These include understandings that have been documented by mathematics education research studies to be essential for learning key ideas of calculus and that reveal whether students actually understand and can use key ideas of precalculus such as average rate of change, function composition, and trigonometric functions. Interview data from CCR validation and other qualitative studies on similar items (Carlson, 1995; Carlson, Oehrtman & Engelke, 2010) revealed that students do not select the correct response on CCR and similar items without having the intended reasoning abilities and understandings. Selection of incorrect responses has also been documented to be associated with weak understandings or misconceptions that have been verified by clinical interviews to be broadly held by precalculus level students and problematic for students in learning calculus.

The CCR is useful as a tool to assess the effectiveness of a precalculus course or curriculum in preparing students for calculus. It can also be used to advise students about their readiness for calculus. We expect that CCR correlations with success in calculus will be higher when administered as a pre-test to students enrolled in calculus courses with a strong conceptual orientation. Even though calculus courses vary in the amount of emphasis placed on skills, techniques, and understanding and using key concepts, we believe that CCR is highly effective in determining the degree to which students are prepared for learning calculus. However, we encourage those who administer CCR to beginning calculus students to use the cut scores that we have reported as advisory, and to consider local constraints and current curriculum foci in precalculus and beginning calculus to adjust break points accordingly.